\RequirePackage{ifpdf}
\ifpdf 
\documentclass[pdftex]{sigma}
\else
\documentclass{sigma}
\fi

\numberwithin{equation}{section}

\begin{document}
\allowdisplaybreaks

\renewcommand{\PaperNumber}{008}

\FirstPageHeading

\renewcommand{\thefootnote}{$\star$}

\ShortArticleName{The Virasoro Algebra and Some Exceptional Lie and Finite Groups}

\ArticleName{The Virasoro Algebra and Some Exceptional Lie\\
 and Finite Groups\footnote{This paper is a contribution to the Proceedings of
the O'Raifeartaigh Symposium on Non-Perturbative and Symmetry
Methods in Field Theory
 (June 22--24, 2006, Budapest, Hungary).
The full collection is available at
\href{http://www.emis.de/journals/SIGMA/LOR2006.html}{http://www.emis.de/journals/SIGMA/LOR2006.html}}}

\Author{Michael P. TUITE}

\AuthorNameForHeading{M.P. Tuite}

\Address{Department of Mathematical Physics, 
National University of Ireland, Galway, Ireland}

\Email{\href{mailto:Michael.Tuite@nuigalway.ie}{Michael.Tuite@nuigalway.ie}}

\ArticleDates{Received October 09, 2006, in f\/inal form December
16, 2006; Published online January 08, 2007}

\Abstract{We describe a number of relationships between properties of the vacuum Verma
module of a Virasoro algebra and the automorphism group of certain vertex
operator algeb\-ras. These groups include the Deligne exceptional series of
simple Lie groups and some exceptional f\/inite simple groups including the
Monster and Baby Monster.}


\Classification{17B68; 20D08; 17B69; 81R05; 81R10}

\section{Introduction}

This paper is based on a talk given at the Lochlainn O'Raifeartaigh
Symposium on Non-Perturbative and Symmetry Methods in Field Theory,
Budapest, 2006. We describe a number of relationships between the vacuum
Verma module of a Virasoro algebra and the Deligne exceptional series of Lie
algebras and also some exceptional f\/inite groups such as the Monster and
Baby Monster. The setting in which this is explained is the theory of Vertex
Operator Algebras (VOAs) or chiral conformal f\/ield theory. In particular, we
construct certain Casimir vectors that are invariant under the VOA
automorphism group. If these vectors are elements of the vacuum Verma module
then the VOA Lie or Griess algebra structure is constrained.  
This is an idea originally introduced by Matsuo for Griess algebras~\cite{M} and
further developed for Lie algebras in~\cite{MMS}. Our approach is more
general with weaker assumptions and is based on a consideration of a certain
rational correlation function. The constraints on the Lie or Griess algebra
structure arise from an analysis of the expansion of this correlation
function in various domains. The paper contains a quite elementary
description of the main ingredients where many of the detailed proofs are
omitted. These will appear elsewhere \cite{T}.

\renewcommand{\thefootnote}{\arabic{footnote}}
\setcounter{footnote}{0}

\section{The Virasoro algebra, Verma modules\\ and the Kac determinant}

We begin with a number of basic def\/initions and properties for Virasoro
algebras e.g.~\cite{KR}. The Virasoro Algebra $\mathrm{Vir}$\ with central
charge $C$ is given by 
\begin{gather}
\lbrack L_{m},L_{n}]=(m-n)L_{m+n}+(m^{3}-m)\frac{C}{12}\delta _{m,-n},\qquad
\lbrack L_{m},C]=0.  \label{Vir}
\end{gather}
The Vacuum Verma Module $V(C,0)$ for $\mathrm{Vir}$\ is def\/ined as follows.
Let $\mathbf{1\in }V(C,0)$ denote the vacuum vector where 
\begin{gather*}
L_{n}\mathbf{1}=0,\qquad n\geq -1.
\end{gather*}
Then $V(C,0)$ is the set of Virasoro descendents of the vacuum (highest
weight) vector 
\begin{gather*}
V(C,0)=\mathbb{C}[L_{-n_{1}}L_{-n_{2}}\cdots L_{-n_{k}}\mathbf{1|}n_{1}\geq n_{2}\geq
\cdots \geq n_{k}\geq 2].
\end{gather*}
$V(C,0)$ is a module for $\mathrm{Vir}$ graded by $L_{0}$ where 
\begin{gather*}
L_{0}(L_{-n_{1}}\cdots L_{-n_{k}}\mathbf{1})=nL_{-n_{1}}\cdots L_{-n_{k}}
\mathbf{1},
\end{gather*}
and where $n=n_{1}+\cdots +n_{k}\geq 0$ is called the (Virasoro) level. Then 
\begin{gather*}
V(C,0)=\bigoplus_{n\geq 0}V^{(n)}(C,0),
\end{gather*}
where $V^{(n)}(C,0)$ denotes the vectors of level $n$. Clearly $\dim
V^{(n)}(C,0)<\infty $ for all $n$.

In general we may consider a Verma Module $V(C,h)$ def\/ined in terms of
highest weight vector $v\in V(C,h)$ obeying 
\begin{gather}
L_{n}v=h\delta _{n,0}v,\qquad n\geq 0.  \label{primary}
\end{gather}%
$v$ is called a primary vector of level $h$. Then $V(C,h)$ is generated by
the Virasoro descendents of $v$. However, we are primarily interested only
in the vacuum module $V(C,0)$ here.

$V(C,0)$ is an irreducible module provided no descendent vector is itself a
primary vector. The irreducibility of $V(C,0)$ can be established by
considering the Kac determinant def\/ined as follows. Def\/ine a symmetric
bilinear form $\langle \cdot,\cdot\rangle $ on $V(C,0)$ with 
\begin{gather*}
\langle L_{-n}u,v\rangle =\langle u,L_{n}v\rangle ,
\end{gather*}%
for $u,v\in V(C,0)$ and normalization $\langle \mathbf{1},\mathbf{1}\rangle
=1$. Note that $\langle u,v\rangle =0$ for $u,v$ of dif\/ferent Virasoro
level. Thus we consider the level $n$ Gram matrix 
\begin{gather*}
M^{(n)}=(\langle u,v\rangle ),
\end{gather*}%
for vacuum descendents $u,v\in V^{(n)}(C,0)$. Then $V(C,0)$ is irreducible
if\/f the Kac determinant $\det M^{(n)}\neq 0$ (with a similar formulation for 
$V(C,h)$) \cite{KR}.

There exists a well-known general product formula for $\det M^{(n)}$ due to
Kac and proved by Feigin and Fuchs \cite{KR}. However, for our purposes, it
is useful to explicitly display $\det M^{(n)}$ for even levels $n\leq 10$.
At level $n=2$ we have $V^{(2)}(C,0)=\mathbb{C}\omega$ where $\omega
=L_{-2}\mathbf{1}$ is called the Conformal Vector. Then 
\begin{gather*}
\det M^{(2)}=\langle \omega ,\omega \rangle =\langle \mathbf{1},L_{2}L_{-2}%
\mathbf{1}\rangle =\frac{C}{2}.
\end{gather*}%
Thus $C\neq 0$ for irreducible $V(C,0)$.

At level $4$ we have $V^{(4)}(C,0)=\mathbb{C}[L_{-2}L_{-2}\mathbf{1},L_{-4}\mathbf{1}]$
with 
\begin{gather*}
M^{(4)} =\left[ 
\begin{array}{cc}
C(4+\frac{1}{2}C) & 3C \\ 
3C & 5C%
\end{array}%
\right] , \qquad
\det M^{(4)} =\frac{1}{2}{C}^{2}\left( 5C+22\right) .
\end{gather*}%
Thus $C\neq 0,-22/5$ for irreducible $V(C,0)$.

Similarly we may compute $\det M^{(n)}$ for $n=6,8,10$ to f\/ind 
\begin{gather*}
\det M^{(6)}= {\frac{3}{4}}\,{C}^{4}\left( 5C+22\right) ^{2}\left(
2\,C-1\right) \left( 7C+68\right) , \\
\det M^{(8)}= {3}\,{C}^{7}\left( 5C+22\right) ^{4}\left( 2C-1\right)
^{2}\left( 7C+68\right) ^{2}\left( 3C+46\right) \left( 5C+3\right) , \\
\det M^{(10)}= {\frac{225}{2}C}^{12}\left( 5C+22\right) ^{8}\left(
2C-1\right) ^{5}\left( 7C+68\right) ^{4}\left( 3C+46\right) ^{2}
 \left( 5C+3\right) ^{2}\left( 11C+232\right) .
\end{gather*}%
where $\dim V^{(6)}(C,0)=4$, $\dim V^{(8)}(C,0)=7$ and $\dim
V^{(10)}(C,0)=12 $.

\section{Some exceptional group numerology}

Let us consider the prime factors of the Kac determinant $\det M^{(n)}$ for
level $n\leq 10$ for \textit{particular values} of $C$. We observe some
coincidences with properties of a number of exceptional Lie and f\/inite
groups. Later on we will explain the underlying reason for these
coincidences and obtain many other relationships with the Virasoro structure.

\textbf{Deligne's exceptional Lie algebras 
$A_{1}$, $A_{2}$, $G_{2}$, $D_{4}$, $F_{4}$, $E_{6}$, $E_{7}$, $E_{8}$.} This set of
simple Lie algebras has been shown recently to share a surprising number of
representation theory properties in common \cite{D}. For example, the
dimension of the adjoint representation $d$ of each of these algebras can be
expressed in terms of the dual Coxeter number $h^{\vee }$ for the given
algebra (as originally found by Vogel): 
\begin{gather}
d=\frac{2(5h^{\vee }-6)(h^{\vee }+1)}{h^{\vee }+6}.  \label{dadj}
\end{gather}

Let us compare $d$ to the independent factors $C$ and $5C+22$ of the level 4
Kac determinant $\det M^{(4)}$ for particular values of $C$:

\begin{center}
{\bf Table 1.}
\vspace{1mm}

\begin{tabular}{|c|c|c|c|c|c|c|c|c|}
\hline
& $A_{1}$ & $A_{2}$ & $G_{2}$ & $D_{4}$ & $F_{4}$ & $E_{6}$ & $E_{7}$ & $%
E_{8}$ \\ \hline
$h^{\vee }$ & $2$ & $3$ & $4$ & $6$ & $9$ & $12$ & $18$ & $30$ \\ \hline
$d$ & $3$ & $2^{3}$ & $2\cdot 7$ & $2^{2}\cdot 7$ & $2^{2}\cdot 13$ & $2\cdot 3\cdot 13$ & $7\cdot 19$ & $%
2^{3}\cdot 31$ \\ \hline
$C$ & $1$ & $2$ & $\frac{2\cdot 7}{5}$ & $2^{2}$ & $\frac{2\cdot 13}{5}$ & $2\cdot 3$ & $7$
& $2^{3}$ \\ \hline
$5C+22$ & $3^{3}$ & $2^{5}$ & $2^{2}\cdot 3^{2}$ & $2\cdot 3\cdot 7$ & $2^{4}\cdot 3$ & $%
2^{2}\cdot 13 $ & $3\cdot 19$ & $2\cdot 31$ \\ \hline
\end{tabular}

\end{center}

The f\/irst row of Table~1 shows $h^{\vee }$ whereas the second row shows the
prime factorization of $d$ for each algebra. The third row shows the prime
factorization for particular values of $C$ (where $C$ is the rank of the
algebra for the simply-laced cases $A_{1}$, $A_{2}$, $D_{4}$, $E_{6}$, $E_{7}$, $E_{8}$).
Notice that each prime divisor of $d$ is a prime divisor of the numerator of
either $C$ or $5C+22$ and hence the numerator of $\det M^{(4)}$.

Some prime divisors of the order of a number of exceptional f\/inite groups
are also related to the Kac determinant factors. We highlight three examples.

\textbf{The Monster simple group $\mathbb{M}$.} The classif\/ication
theorem of f\/inite simple groups states that a f\/inite simple group is either
one of several inf\/inite families of simple groups (e.g.~the alternating
groups $A_{n}$ for $n\geq 5$) or else is one of 26 sporadic f\/inite simple
groups. The largest sporadic group is the Monster group $\mathbb{M}$ of
(prime factored) order 
\begin{gather*}
|\mathbb{M}%
|=2^{46}\cdot 3^{20}\cdot 5^{9}\cdot 7^{6}\cdot 11^{2}\cdot 13^{3}\cdot 17\cdot 19\cdot 23\cdot
29\cdot 31\cdot 41\cdot 47\cdot 59\cdot 71.
\end{gather*}%
The two lowest dimensional non-trivial irreducible representations are of
dimension $d_{2}=196883$ $=47\cdot 59\cdot 71$ and $d_{3}=21296876=2^{2}\cdot 31\cdot 41\cdot 59\cdot 71$.
Consider the independent factors of $\det M^{(10)}$ for $C=24$

\begin{center}
{\bf Table 2.}

\vspace{1mm}

$ 
\begin{tabular}{|c|c|c|c|c|c|c|}
\hline
${C}$ & $5C+22$ & $2C-1$ & $7C+68$ & $3C+46$ & $5C+3$ & $11C+232$ \\ \hline
$2^{3}\cdot 3$ & $2\cdot 71$ & $47$ & $2^{2}\cdot 59$ & $2\cdot 59$ & $3\cdot 41$ & $2^{4}\cdot 31$ \\ 
\hline
\end{tabular}
$

\end{center}

Notice that the prime divisors $2$, $31$, $41$, $47$, $59$, $71$ of $|\mathbb{M}|$ (which
are all of the prime divisors of $d_{2}$ and $d_{3}$) are divisors of $\det
M^{(10)}$.

\textbf{The Baby Monster simple group $\mathbb{B}$.} The second
largest sporadic group is the Baby Monster group $\mathbb{B}$ of order 
\begin{gather*}
|\mathbb{B}|=2^{41}\cdot 3^{13}\cdot 5^{6}\cdot 7^{2}\cdot 11\cdot 13\cdot 17\cdot 19\cdot 23\cdot 31\cdot 47.
\end{gather*}

The prime divisors $2$, $3$, $5$, $23$, $31$, $47$ of $|\mathbb{B}|$ are divisors of the
numerator of the independent factors of $\det M^{(6)}$ for $C=23\frac{1}{2}$:

\begin{center}
{\bf Table 3.}
\vspace{1mm}

$ 
\begin{tabular}{|c|c|c|c|}
\hline
${C}$ & $5C+22$ & $2C-1$ & $7C+68$ \\ \hline
$\frac{47}{2}$ & $\frac{3^{2}\cdot 31}{2}$ & $2\cdot 23$ & $\frac{3\cdot 5\cdot 31}{2}$ \\ \hline
\end{tabular}
$

\end{center}

\textbf{The simple group $O_{10}^{+}(2)$.} This group has order 
\begin{gather*}
|O_{10}^{+}(2)|=2^{20}\cdot 3^{5}\cdot 5^{2}\cdot 7\cdot 17\cdot 31.
\end{gather*}

The prime divisors $2$, $3$, $5$, $31$ of $|O_{10}^{+}(2)|$ are divisors of the
independent factors of $\det M^{(6)}$ for $C=8$

\begin{center}
{\bf Table 4.}
\vspace{1mm}

$ 
\begin{tabular}{|c|c|c|c|}
\hline
${C}$ & $5C+22$ & $2C-1$ & $7C+68$ \\ \hline
$2^{3}$ & $2\cdot 31$ & $3\cdot 5$ & $2^{2}\cdot 31$ \\ \hline
\end{tabular}
$
\end{center}

\section{Vertex operator algebras}

\subsection{Axioms}

The observations made in the last section can be understood in the context
of Vertex Operator Algebras (VOAs) \cite{B,Ka,MN}. The basic
idea is that the groups appearing above arise as automorphism (sub)groups of
particular VOAs with the given central charge $C$. In this section we review
the relevant aspects of VOA theory required in the subsequent sections. 

A Vertex Operator Algebra (VOA) consists of a $\mathbb{Z}$-graded vector
space $V=\bigoplus_{k\geq 0}V^{(k)}$ with $\dim V^{(k)}<\infty $ and with
the following properties:

\begin{itemize}
\item \textbf{Vacuum.} $V^{(0)}=\mathbb{C}\mathbf{1}$ for vacuum vector $\mathbf{1}
$.

\item \textbf{Vertex operators (state-f\/ield correspondence).} For each $a\in
V^{(k)}$ we have a~vertex operator 
\begin{gather}
Y(a,z)=\sum_{n\in \mathbb{Z}}a_{n}z^{-n-k},  \label{Y(a,z)}
\end{gather}
with component operators (modes) $a_{n}\in \mathrm{End}V$ such that 
\begin{gather*}
Y(a,z)\cdot \mathbf{1|}_{z=0}=a_{-k}\cdot \mathbf{1}=a.  
\end{gather*}
$z$ is a formal variable here (but is taken as a complex number in physics).
Note also that we are employing ``physics modes'' in (\ref{Y(a,z)}).

The vacuum vector has vertex operator 
\begin{gather*}
Y(\mathbf{1},z)=Id_{V},  
\end{gather*}%
so that $\mathbf{1}_{n}=\delta _{n,0}$.

\item \textbf{Virasoro subalgebra.} There exists a conformal vector $\omega
\in V^{(2)}$ with 
\begin{gather*}
Y(\omega ,z)=\sum_{n\in \mathbb{Z}}L_{n}z^{-n-2},  
\end{gather*}%
where the modes $L_{n}$ form a Virasoro algebra (\ref{Vir}) of central
charge $C$. The $\mathbb{Z}$-grading is determined by $L_{0}$ with $%
V^{(k)}=\mathbb{C}[a\in V|L_{0}a=ka]$.

$L_{-1}$\textbf{\ }acts as a translation operator with 
\begin{gather*}
Y(L_{-1}a,z)=\partial _{z}Y(a,z).  
\end{gather*}

\item \textbf{Locality.} For any pair of vertex operators we have for
suf\/f\/iciently large integer $N$ 
\begin{gather}
(x-y)^{N}[Y(a,x),Y(b,y)]=0.  \label{Locality}
\end{gather}
\end{itemize}

These axioms lead to the following basic VOA properties: e.g.~\cite{Ka,MN}.

\textbf{Translation.} For any $a\in V$ then for $\left\vert y\right\vert
<\left\vert x\right\vert $ (i.e. formally expanding in $y/x$) 
\begin{gather}
e^{yL_{-1}}Y(a,x)e^{-yL_{-1}}=Y(a,x+y).  \label{Translation}
\end{gather}

\textbf{Skew-symmetry. }For $a,b\in V$ then 
\begin{gather}
Y(a,z)b=e^{zL_{-1}}Y(b,-z)a.  \label{Skew-sym}
\end{gather}

\textbf{Associativity.} For $a,b\in V$ then for $\left\vert x-y\right\vert
<\left\vert y\right\vert <\left\vert x\right\vert $ we f\/ind\footnote{Strictly speaking, (\ref{Associativity})
holds for matrix elements taken with repect to any pair $u'$, $v$ for any $v\in V$ and $u'\in V'$, the dual
vector space~\cite{FHL}.}
\begin{gather}
Y(a,x)Y(b,y)=Y(Y(a,x-y)b,y).  \label{Associativity}
\end{gather}

\textbf{Borcherds' commutator formula.} For $a\in V^{(k)}$ and $b\in V$ then 
\begin{gather}
\lbrack a_{m},b_{n}]=\sum\limits_{j\geq 0}\binom{m+k-1}{j}(a_{j-k+1}b)_{m+n}.
\label{Borcherds}
\end{gather}%
\textbf{Example.} For $a=\omega \in V^{(2)}$ and $m=0$ we have $%
[L_{0},b_{n}]=(L_{-1}b)_{n}+(L_{0}b)_{n}=-nb_{n}$. Thus 
\begin{gather}
b_{n}:\ V^{(m)}\rightarrow V^{(m-n)}.  \label{bnVm}
\end{gather}%
In particular, the zero mode $b_{0}$ is a linear operator on $V^{(m)}$.
Similarly, for all $m$ and any primary vector $b\in V^{(h)}$, we have the
familiar property from conformal f\/ield theory 
\begin{gather}
\lbrack L_{m},b_{n}]=((h-1)m-n)b_{m+n}.  \label{Lmbn}
\end{gather}

\subsection[The Li-Zamolodchikov metric]{The Li--Zamolodchikov metric}

Assume that $V^{(0)}=\mathbb{C}\mathbf{1}$ and $L_{1}v=0$ for all $v\in V^{(1)}$.
Then there exists a unique invariant bilinear form $\langle \cdot ,\cdot \rangle $ where
for all $a,b,c\in V$  \cite{FHL,L} 
\begin{gather*}
\langle Y\left(e^{zL_{1}}\left(-\frac{1}{z^{2}}\right)^{L_{0}}c,\frac{1}{z}\right)a,b\rangle
=\langle a,Y(c,z)b\rangle ,  
\end{gather*}%
and with normalization $\langle \mathbf{1},\mathbf{1}\rangle =1$. For any
quasi-primary $c\in V^{(k)}$ (i.e.\ where $L_{1}c=0$) then 
\begin{gather}
\langle c_{n}a,b\rangle =(-1)^{k}\langle a,c_{-n}b\rangle .  \label{cnab}
\end{gather}%
Choosing $a\in V^{(k)}$, $n=k$ and $b=\mathbf{1}$, (\ref{cnab}) implies for
quasi-primary $c$ 
\begin{gather}
c_{k}a=(-1)^{k}\langle a,c\rangle \mathbf{1}.  \label{cka}
\end{gather}%
For $c=\omega $, (\ref{cnab}) also implies 
\begin{gather*}
\langle L_{n}a,b\rangle =\langle a,L_{-n}b\rangle .  
\end{gather*}%
$\langle \cdot ,\cdot \rangle $ is symmetric \cite{FHL} and, furthermore, is
non-degenerate if\/f $V$ is semisimple~\cite{L}. We call such a~unique
non-degenerate form the Li--Zamolodchikov or Li-Z metric. In particular,
consi\-de\-ring the Li-Z metric on the vacuum Virasoro descendents $%
V(C,0)\subset V$ implies that the Kac determinant $\det M^{(n)}\neq 0$ for
each level $n$.

\subsection{Lie and Griess algebras}

\textbf{Lie algebras.} Let us consider a number of relevant subalgebras of a
VOA $V$ with Li-Z metric. Suppose that $\dim V^{(1)}>0$ and def\/ine for all $%
a,b\in V^{(1)}$ 
\begin{gather}
{\rm ad}\,(a)b=a_{0}b=-b_{0}a,  \label{adab}
\end{gather}%
where the second equality follows from skew-symmetry (\ref{Skew-sym}). We
may thus def\/ine a bracket $[a,b]={\rm ad}\,(a)b$ which satisf\/ies the Jacobi identity
so that $V^{(1)}$ forms a Lie algebra. Furthermore, the Li-Z metric $\langle
a,b\rangle $ on $V^{(1)}$ is an invariant non-degenerate symmetric bilinear
form with\footnote{%
Non-degeneracy of the metric is not necessary here.} 
\begin{gather}
\langle \lbrack a,b],c\rangle =\langle -b_{0}a,c\rangle =\langle
a,b_{0}c\rangle =\langle a,[b,c]\rangle ,  \label{Lie form}
\end{gather}%
for all \ $a,b,c\in V^{(1)}$. The Borcherds' commutator formula (\ref%
{Borcherds}) then def\/ines a Kac--Moody 
algebra\footnote{The non-standard minus sign on the RHS of (\ref{KacMoody}) is due to the
normalization $\langle \mathbf{1},\mathbf{1}\rangle =1$.} 
\begin{gather}
\lbrack a_{m},b_{n}]=(a_{0}b)_{m+n}+(a_{1}b)_{m+n}=[a,b]_{m+n}-m\langle
a,b\rangle \delta _{m+n,0},  \label{KacMoody}
\end{gather}%
using $a_{1}b=-\langle a,b\rangle \mathbf{1}$ from (\ref{cka}).

\textbf{Griess algebras.} Suppose that $\dim V^{(1)}=0$ and consider $a,b\in
V^{(2)}$. Skew-symmetry implies $a_{0}b=b_{0}a$ so that we may def\/ine 
\begin{gather*}
a\bullet b=b\bullet a=a_{0}b, 
\end{gather*}%
to form a commutative non-associative algebra on $V^{(2)}$ known as a Griess
algebra. The Li-Z metric $\langle a,b\rangle $ on $V^{(2)}$ is an invariant
bilinear form with 
\begin{gather*}
\langle a\bullet b,c\rangle =\langle b,a\bullet c\rangle ,
\end{gather*}%
for all\ $a,b,c\in V^{(2)}$.

\subsection{The automorphism group of a VOA}

The automorphism group $\textrm{Aut}\,(V)$ of a VOA $V$ consists of all
elements $g\in GL(V)$ for which 
\begin{gather*}
gY(a,z)g^{-1}=Y(ga,z),  
\end{gather*}%
with $g\omega =\omega $, the conformal vector. Thus the $L_{0}$ grading is
preserved by $\textrm{Aut}\,(V)$ and every Virasoro descendent of the vacuum is
invariant under $\textrm{Aut}\,(V)$. Furthermore, the Li-Z metric is invariant with 
\begin{gather}
\langle ga,gb\rangle =\langle a,b\rangle ,  \label{LiZinv}
\end{gather}%
for all $a,b\in V$ and $g\in \textrm{Aut}\,(V)$.

For VOAs with $\dim V^{(1)}>0$, then $\textrm{Aut}\,(V)$ contains continuous symmetries 
$g=\exp (a_{0})$ generated by elements of the Lie algebra $a\in V^{(1)}$.
VOAs with $\dim V^{(1)}=0$ for which $V^{(2)}$ def\/ines a Griess algebra are
particularly interesting. Examples include the Moonshine Module~$V^{\natural}$ 
of central charge $C=24$ with $\textrm{Aut}\,(V)=\mathbb{M}$, the Monster
group \cite{FLM}. In this case, $V^{\natural (2)}$ is the original Griess
algebra of dimension $196884=1+196883$ where $V^{\natural (2)}$ decomposes
into the $\mathbb{M}$ invariant Virasoro vector $\omega $ and an irreducible 
$\mathbb{M}{\ }$representation of dimension $196883=47\cdot 59\cdot 71$. Other
examples include a VOA with $C=23\frac{1}{2}$ and $\textrm{Aut}\,(V)=\mathbb{B}
$, the Baby Monster where $\dim V^{(2)}=1+96255$ with an irreducible $%
\mathbb{B}$ representation of dimension $96255=3^{3}\cdot 5\cdot 23\cdot 31$~\cite{H} and
a VOA with $C=8$ and \textrm{Aut}$(V)=O_{10}^{+}(2)\cdot 2$ and $\dim V^{(2)}=1+155$ 
with an irreducible $O_{10}^{+}(2)$
representation of dimension $155=5\cdot 31$~\cite{G}.  In each case, $C$
corresponds to the values shown in Tables~2, 3 and 4.

\section{Lie algebras and Virasoro descendents}

\subsection{Quadratic Casimirs}

We now discuss the relationship between the structure of the Lie algebra $%
V^{(1)}$ and the Virasoro algebra described earlier. We initially follow a
technique due to Matsuo \cite{M,MMS} for constructing a Casimir
vector $\lambda ^{(n)}\in V^{(n)}$ which is $\textrm{Aut}\,(V)$ invariant. The
relationship with the Virasoro algebra follow provided $\lambda ^{(n)}\in
V^{(n)}(C,0)$, the vacuum Virasoro descendents of level $n$. This occurs,
for example, if $V^{(n)}$ contains no $\mathrm{Aut}\,(V)$ invariants apart
from the elements of $V^{(n)}(C,0)$. Here we describe a new general method
based on weaker assumptions for obtaining various constraints on $V^{(1)}$
by considering the expansion in various domains of a particular correlation
function~\cite{T}.

Consider a VOA with Li-Z metric with 
\begin{gather*}
d=\dim V^{(1)}>0,  
\end{gather*}%
the dimension of the Lie algebra $V^{(1)}$. Let $V^{(1)}$ have basis 
$\{u_{\alpha }\,|\,\alpha =1,\ldots, d\}$ and dual basis $\{u^{\beta }|\beta
=1,\ldots, d\}$ i.e.\ $\langle u^{\alpha },u_{\beta }\rangle =\delta _{\beta
}^{\alpha }$. For $n\geq 0$ def\/ine the quadratic Casimir vector \cite{M} 
\begin{gather*}
\lambda ^{(n)}=u_{1-n}^{\alpha }u_{\alpha }\in V^{(n)},  
\end{gather*}%
where $\alpha $ is summed. $\lambda ^{(n)}$ is invariant under $\textrm{Aut}\,
(V)$ following (\ref{LiZinv}). We f\/ind from (\ref{cka}) and~(\ref{adab})
that 
\begin{gather*}
\lambda ^{(0)} =-d\mathbf{1},\qquad 
\lambda ^{(1)} =0. 
\end{gather*}%
Furthermore, using (\ref{Lmbn}) it follows that for all $m>0$
\begin{gather}
L_{m}\lambda ^{(n)}=(n-1)\lambda ^{(n-m)}.  \label{Lmlambda}
\end{gather}

Suppose that $\lambda ^{(n)}\in V^{(n)}(C,0)$ i.e.\ $\lambda ^{(n)}$ is a
Virasoro descendent of the vacuum of level~$n$. Then $\lambda ^{(n)}
$ can be determined exactly via~(\ref{Lmlambda}) \cite{M}. For
example, if $\lambda ^{(2)}\in V^{(2)}(C,0)$ then $\lambda ^{(2)}=k\omega $
for some $k$. Hence 
\begin{gather*}
\langle \omega ,\lambda ^{(2)}\rangle =k\langle \omega ,\omega \rangle =k%
\frac{C}{2}.
\end{gather*}%
But from (\ref{cka}) and (\ref{Lmlambda}) we f\/ind $\langle \omega ,\lambda
^{(2)}\rangle \mathbf{1}=L_{2}\lambda ^{(2)}=\lambda ^{(0)}=-d\mathbf{1}$.
Thus, for $C\neq 0$ 
\begin{gather}
\lambda ^{(2)}=-\frac{2d}{C}\omega .  \label{lambda2}
\end{gather}%
$\lambda ^{(2)}$ is necessarily singular at $C=0$ where $\det M^{(2)}=0$.
Notice that (\ref{lambda2}) is just the  familiar Sugawara construction of $%
\omega $ in the case that $V^{(1)}$ is a simple Lie algebra.

Assuming $\lambda ^{(4)}\in V^{(4)}(C,0)$ we similarly f\/ind 
\begin{gather}
\lambda ^{(4)}=\frac{3d}{C\left( 5C+22\right) }[4L_{-2}L_{-2}\mathbf{1}%
+\left( 2+C\right) L_{-4}\mathbf{1}],  \label{lambda4}
\end{gather}%
which is singular at $C=0,-22/5$ where $\det M^{(4)}=0$. 

\subsection[Rational correlation functions and the $V^{(1)}$ Killing form]{Rational 
correlation functions and the $\boldsymbol{V^{(1)}}$ Killing form}

Let $a,b \in V^{(1)}$ and consider the following correlation function 
\begin{gather*}
F(a,b;x,y)=\langle a,Y(u^{\alpha },x)Y(u_{\alpha },y)b\rangle .  
\end{gather*}%
Then we have

\begin{proposition}
\label{Proposition_F_Lie} $F(a,b;x,y)$ is a rational function 
\begin{gather}
F(a,b;x,y)=\frac{G(a,b;x,y)}{x^{2}y^{2}(x-y)^{2}},  \label{Fxy_rat}
\end{gather}%
where $G(a,b;x,y)$ is bilinear in $a$, $b$ and is a homogeneous, symmetric
polynomial in $x$, $y$ of degree $4$.
\end{proposition}

\begin{proof} The bilinearity of $G$ in $a$, $b$ is obvious. Locality 
(\ref{Locality}) implies that $F(a,b;x,y)$ is of the form (\ref{Fxy_rat}) \cite{FHL,MN} 
where $G(a,b;x,y)$ is clearly symmetric. The degree and
homogeneity follow from (\ref{Y(a,z)}) and (\ref{bnVm}).
\end{proof}

It is convenient to parameterize $G(a,b;x,y)$ in terms of $3$ independent
coef\/f\/icients $P(a,b)$, $Q(a,b)$, $R(a,b)$ as follows 
\begin{gather*}
G(a,b;x,y)=P(a,b)x^{2}y^{2}+Q(a,b)xy(x-y)^{2}+R(a,b)(x-y)^{4}.  
\end{gather*}%
Expanding the rational function $F(a,b;x,y))$ in $\frac{x-y}{y}$ we f\/ind 
\begin{gather}
F(a,b;x,y)=(x-y)^{-2}\left( P(a,b)+Q(a,b)\left( \frac{x-y}{y}\right)
^{2}+\cdots \right) ,  \label{Fab_gexpan1}
\end{gather}%
whereas expanding in $\frac{-y}{x-y}$ we f\/ind 
\begin{gather}
F(a,b;x,y)=y^{-2}\left( R(a,b)+(2R(a,b)-Q(a,b))\left( \frac{-y}{x-y}\right)
+\cdots \right) .  \label{Fab_gexpan2}
\end{gather}%
From the corresponding VOA expansions, $P$, $Q$, $R$ may be computed in terms of
the Li-Z metric $\langle a,b\rangle $ and the Lie algebra Killing form 
\begin{gather*}
K(a,b)={\rm Tr}_{V^{(1)}}({\rm ad}(a){\rm ad}(b))={\rm Tr}_{V^{(1)}}(a_{0}b_{0}),
\end{gather*}%
as follows:

\begin{proposition}
\label{Proposition_Lie} $P(a,b)$, $Q(a,b)$, $R(a,b)$ are given by 
\begin{gather*}
P(a,b) =-d\langle a,b\rangle ,  \\
Q(a,b) =K(a,b)-2\langle a,b\rangle ,  \\
R(a,b) =-\langle a,b\rangle .  
\end{gather*}
\end{proposition}

\begin{proof} 
Associativity (\ref{Associativity}) implies we may expand in 
$(x-y)/y$ to obtain%
\begin{gather}
F(a,b;x,y) =\langle a,Y(Y(u^{\alpha },x-y)u_{\alpha },y)b\rangle  
=\sum_{n\geq 0}\langle a,Y(\lambda ^{(n)},y)b\rangle (x-y)^{n-2} 
\notag \\
\phantom{F(a,b;x,y)}{}=(x-y)^{-2}\sum_{n\geq 0}\langle a,\lambda _{0}^{(n)}b\rangle
\left( \frac{x-y}{y}\right) ^{n}.  \label{Fxy_expan1}
\end{gather}%
The leading term of (\ref{Fxy_expan1}) is determined by $\lambda
_{0}^{(0)}=-d$ so that comparing with (\ref{Fab_gexpan1}) we f\/ind $%
P(a,b)=-d\langle a,b\rangle $. Note that the next to leading term
automatically vanishes since $\lambda _{0}^{(1)}=0$.

We may alternatively expand $F(a,b;x,y)$ in $-y/(x-y)$ to f\/ind 
\begin{gather}
F(a,b;x,y) =\langle a,Y(u^{\alpha },x)e^{yL_{-1}}Y(b,-y)u_{\alpha }\rangle 
\notag \\
\phantom{F(a,b;x,y)}{}=\langle a,e^{yL_{-1}}Y(u^{\alpha },x-y)Y(b,-y)u_{\alpha }\rangle   \notag
\\
\phantom{F(a,b;x,y)}{}=\langle e^{yL_{1}}a,Y(u^{\alpha },x-y)Y(b,-y)u_{\alpha }\rangle   \notag
\\
\phantom{F(a,b;x,y)}{}=\langle a,Y(u^{\alpha },x-y)Y(b,-y)u_{\alpha }\rangle   \notag \\
\phantom{F(a,b;x,y)}{}=y^{-2}\sum_{m\geq 0}\langle a,u_{m-1}^{\alpha }b_{1-m}u_{\alpha }\rangle
\left( \frac{-y}{x-y}\right) ^{m},  \label{Fxy_expan2}
\end{gather}%
(respectively using skew-symmetry (\ref{Skew-sym}), translation (\ref{Translation}), 
invariance of the Li-Z metric (\ref{cnab}) and that $a$ is
primary (\ref{primary})). Comparing to the leading term of (\ref{Fab_gexpan2}) we f\/ind using (\ref{cka}) that 
\begin{gather*}
R(a,b) =\langle a,u_{-1}^{\alpha }b_{1}u_{\alpha }\rangle  
=-\langle a,u^{\alpha }\rangle \langle b,u_{\alpha }\rangle =-\langle
a,b\rangle .
\end{gather*}%
Comparing the next to leading terms of (\ref{Fab_gexpan2}) and (\ref%
{Fxy_expan2}) and using (\ref{adab}) and (\ref{Lie form}) we f\/ind 
\begin{gather*}
2R(a,b)-Q(a,b) =\langle a,u_{0}^{\alpha }b_{0}u_{\alpha }\rangle 
=-\langle u^{\alpha },a_{0}b_{0}u_{\alpha }\rangle =-K(a,b). \tag*{\qed}
\end{gather*}
\renewcommand{\qed}{}
\end{proof}

We next show that if $\lambda ^{(2)}\,$is a vacuum Virasoro descendent then
the Killing form is proportional to the Li-Z metric:

\begin{proposition}
\label{Proposition_Lie_lambda2} Suppose that $\lambda ^{(2)}\in V^{(2)}(C,0)$. Then 
\begin{gather}
K(a,b)=-2\langle a,b\rangle \left( \frac{d}{C}-1\right).
\label{KillingForm}
\end{gather}
\end{proposition}

\begin{proof}
If $\lambda ^{(2)}\in V^{(2)}(C,0)$ then (\ref{lambda2})
implies $\lambda _{0}^{(2)}=-\frac{2d}{C}L_{0}$. The $n=2$ term in (\ref%
{Fxy_expan1}) is thus $-\frac{2d}{C}\langle a,b\rangle $. Comparing to (\ref%
{Fab_gexpan1}) we f\/ind 
\begin{gather*}
Q(a,b)=-\frac{2d}{C}\langle a,b\rangle ,
\end{gather*}%
and hence the result follows from Proposition \ref{Proposition_Lie}. 
\end{proof}

Since the Li-Z metric is non-degenerate, it immediately follows from
Cartan's condition that

\begin{corollary}
\label{Corollary_Killing} The Lie algebra $V^{(1)}$ is semi-simple for $%
d\neq C$ and is solvable for $d=C$.
\end{corollary}

The Killing form (\ref{KillingForm}) has previously arisen in the literature
from considerations of modu\-lar invariance in the classif\/ication of $V^{(1)}$
for holomorphic self-dual VOAs in the work of Schellekens~\cite{S} for
central charge $C=24$ and Dong and Mason \cite{DM} for $C=8,16,24$. A~similar 
result also appears in~\cite{MMS} also based on a Casimir
invariant approach.

We can repeat results of~\cite{S,DM} concerning the
decomposition of $V^{(1)}$ into simple components for $d\neq C$. Corollary %
\ref{Corollary_Killing} implies 
\begin{gather}
V^{(1)}=\mathfrak{g}_{1,k_{1}}\oplus \mathfrak{g}_{2,k_{2}}\oplus \mathfrak{%
\cdots }\oplus \mathfrak{g}_{r,k_{r}},  \label{V1simple}
\end{gather}%
where $\mathfrak{g}_{i,k_{i}}$ is a simple Kac--Moody algebra of dimension $%
d_{i}$ (where $d=\sum\limits_{i=1,\dots,r}d_{i}$) and level 
\begin{gather*}
k_{i}=-\frac{1}{2}\langle \mathbf{\alpha }_{i},\mathbf{\alpha }_{i}\rangle ,
\end{gather*}%
where $\mathbf{\alpha }_{i}$ is a long root of $\mathfrak{g}_{i,k_{i}}$ (so
that $(a,b)_{i}=-\langle a,b\rangle /k_{i}$ def\/ines a non-degenerate form 
on~$\mathfrak{g}_{i,k_{i}}$ with normalization $(\mathbf{\alpha }_{i},\mathbf{%
\alpha }_{i})_{i}=2$). The dual Coxeter number $h_{i}^{\vee }$ for $%
\mathfrak{g}_{i,k_{i}}$ is then found from the Killing form $K(h_{\mathbf{%
\alpha }_{i}},h_{\mathbf{\alpha }_{i}})=4h_{i}^{\vee }$. Hence (\ref%
{KillingForm}) implies that for each simple component 
\begin{gather*}
\frac{h_{i}^{\vee }}{k_{i}}=\frac{d}{C}-1.  
\end{gather*}%
This implies that a f\/inite number of solutions for (\ref{V1simple}) exist
for any pair $(C,d)$.

\subsection{Deligne's exceptional Lie algebras}

We next consider the expansion (\ref{Fxy_expan2}) of $F(a,b;x,y)$ to the
next leading term assuming that $\lambda ^{(4)}$ is also a vacuum Virasoro
descendent so that (\ref{lambda4}) holds. This results in one further
constraint on $P$, $Q$, $R$ leading to~\cite{T}

\begin{proposition}
\label{Proposition_d(c)_Lie} Suppose that $\lambda ^{(n)}\in V^{(n)}(C,0)$
for $n\leq 4$. Then 
\begin{gather*}
d(C)=\frac{C\left( 5C+22\right) }{10-C}.  
\end{gather*}
\end{proposition}

Note that $d(C)$necessarily contains the independent factors of the level 
$4$ Kac determinant $\det M^{(4)}$ since it must vanish for $C=0,-22/5$ for
which the construction of $\lambda ^{(4)}\,$is impossible.

For positive rational central charge $C$ then $d(C)$ is a positive integer
for only 21 dif\/ferent values of $C$. The remarks subsequent to Corollary \ref%
{Corollary_Killing} imply that the possible values are further restricted
and that $V^{(1)}$ must be one of the Deligne exceptional simple Lie
algebras \cite{T}:

\begin{proposition}
\label{Proposition_Deligne}Suppose that $\lambda ^{(n)}\in V^{(n)}(C,0)$ for 
$n\leq 4$ and that $C\in \mathbb{Q}_{+}\,$. Then $V^{(1)}$ is one of the
simple Lie algebras $A_{1}$, $A_{2}$, $G_{2}$, $D_{4}$, $F_{4}$, $E_{6}$, $E_{7}$, $E_{8}$ for $%
C=1,2,\frac{14}{5},4,\frac{26}{5},6,7,8$ respectively, with Kac--Moody
level $k=1$ and dual Coxeter number 
\begin{gather}
h^{\vee }(C)=6{\frac{2+C}{10-C}=\frac{d(C)}{C}-1.}
\label{hhat(C)}
\end{gather}
\end{proposition}

$d(C)$ is precisely the original Vogel formula (\ref{dadj}) with dual
Coxeter number given by (\ref{hhat(C)}). The properties of Table 1 are now
obvious since if a prime $p$ divides $d(C)$ then $p$ must divide one of the
factors $C$ or $22+5C$ of $\det M^{(4)}$. We note that results similar to
Propositions \ref{Proposition_d(c)_Lie} and~\ref{Proposition_Deligne} also
appear in~\cite{MMS} but are based on a number of further technical
assumptions.

We conclude this section with a result concerning the constraints arising
from higher level Casimirs being vacuum descendents \cite{T}:

\begin{proposition}
If $\lambda ^{(n)}\in V^{(n)}(C,0)$ for $n\leq 6$ then $V^{(1)}$ must be
either $A_{1}$ or $E_{8} $ with level~$1$. If $\lambda ^{(n)}\in V^{(n)}(C,0)$
for all $n$ then $V^{(1)}$ is $A_{1} $ level $1$.
\end{proposition}

\section{Griess algebras and Virasoro descendents}

We can now repeat the approach taken in the last section for a VOA with a
Li-Z metric with $\dim V^{(1)}=0$ so that $V^{(2)}$ def\/ines a Griess
algebra. This was the original case considered by Matsuo~\cite{M} but was
based on stronger assumptions (e.g.\ assuming the existence of a proper
idempotent or that $\mathrm{Aut}\,(V)$ is f\/inite). We give a brief description
of the constraints arising from quadratic Casimirs being Virasoro
descendents using a similar approach to the last section. Detailed proofs
will appear elsewhere~\cite{T}.

Let$\ \hat{V}^{(2)}=V^{(2)}-\mathbb{C}\omega$ denote the level $2$ primary
vectors of dimension 
\begin{gather*}
d_{2}=\dim \hat{V}^{(2)}>0,
\end{gather*}
with basis $\{u_{\alpha }\}$ and dual basis $\{u^{\alpha }\}$. We def\/ine 
\textrm{Aut}$(V)$ invariant quadratic Casimir vectors 
\begin{gather*}
\mu ^{(n)}=u_{2-n}^{\alpha }u_{\alpha }\in V^{(n)}.  
\end{gather*}
Then 
\begin{gather*}
\mu ^{(0)}=d_{2}\mathbf{1},\qquad \mu ^{(1)}=0,  
\end{gather*}
and for $m>0$%
\begin{gather*}
L_{m}\mu ^{(n)}=(m+n-2)\mu ^{(n-m)}.  
\end{gather*}

Consider the correlation function for $a,b\in \hat{V}^{(2)}$%
\begin{gather*}
F(a,b;x,y)=\langle a,Y(u^{\alpha },x)Y(u_{\alpha },y)b\rangle .
\end{gather*}%
Then we have

\begin{proposition}
$F(a,b;x,y)$ is a rational function 
\begin{gather*}
F(a,b;x,y)=\frac{G(a,b;x,y)}{x^{4}y^{4}(x-y)^{4}},
\end{gather*}%
where $G(a,b;x,y)$ is bilinear in $a$, $b$ and is a homogeneous, symmetric
polynomial in $x$, $y$ of degree~$8$.
\end{proposition}

In this case $G(a,b;x,y)$ is determined by $5$ independent coef\/f\/icients. We
use associativity~(\ref{Associativity}) as in~(\ref{Fxy_expan1}) to expand
in $(x-y)/y$ so that 
\begin{gather*}
F(a,b;x,y)=(x-y)^{-4}\sum\nolimits_{n\geq 0}\langle a,\mu _{0}^{(n)}b\rangle
\left( \frac{x-y}{y}\right) ^{n}.
\end{gather*}%
Similarly to (\ref{Fxy_expan2}), we expand in $-y/(x-y)$ to obtain

\begin{gather*}
F(a,b;x,y)=y^{-4}\left( \langle a,b\rangle +0+Tr_{\hat{V}^{(2)}}(a_{0}b_{0})%
\left( \frac{-y}{x-y}\right) ^{2}+\cdots \right).
\end{gather*}%
Analogously to Propositions \ref{Proposition_Lie} and \ref{Proposition_Lie_lambda2} 
and employing methods of \cite{DM} we f\/ind

\begin{proposition}
Suppose that $\mu ^{(n)}\in V^{(n)}(C,0)$ for $n\leq 4$. Then

{\samepage
\begin{enumerate}\itemsep=0pt
\item[{\rm 1.}] $G(a,b;x,y)$ is given by 
\begin{gather*}
 \langle a,b\rangle \Bigg[ d_{2}{x}^{4}{y}^{4}+\frac{8d_{2}}{C}{x}^{3}{y}%
^{3}{(x-y)}^{2}+\frac{4d_{2}\left( 44-C\right) }{C\left( 5C+22\right) }{x}%
^{2}{y}^{2}{(x-y)}^{4} \\
\qquad{} +2({x}^{2}+{y}^{2}){(x-y)}^{6}-{(x-y)}^{8}\Bigg].
\end{gather*}

\item[{\rm 2.}] $V^{(2)}$ has a non-degenerate trace form 
\begin{gather*}
{\rm Tr}_{V^{(2)}}((a\bullet b)_{0})=\frac{8(d_{2}+1)}{C}\langle a,b\rangle .
\end{gather*}

\item[{\rm 3.}] The Griess algebra $V^{(2)}$ is simple.
\end{enumerate}}
\end{proposition}
 
$d_{2}$ and $C$ are related by the next Virasoro constraint in analogy to
Proposition \ref{Proposition_d(c)_Lie}:

\begin{proposition}
Suppose that $\mu ^{(n)}\in V^{(n)}(C,0)$ for $n\leq 6$. Then 
\begin{gather*}
d_{2}(C)=\frac{1}{2}\frac{\left( 5C+22\right) \left( 2C-1\right) \left(
7C+68\right) }{748-55C+C^{2}}.  
\end{gather*}
\end{proposition}

The formula for $d_{2}(C)+1$ previously appeared in~\cite{M} subject to
further assumptions. The appearance of the independent factors of $\det
M^{(6)}$ in $d_{2}(C)$ follows as before. (The absence of a $C$ factor
follows from the assumed decomposition $V^{(2)}=\mathbb{C}\omega\oplus \hat{V}%
^{(2)}$). If a prime $p$ divides $d_{2}(C)$ then $p$ must divide (at least)
one of the independent factors $5C+22$, $2C-1$ or $7C+68$ of $\det M^{(6)}$.
This explains many of the prime divisor properties (but not all) for Tables
2, 3 and~4.

There are exactly 36 positive rational values for $C$ for which $d_{2}(C)\,$%
is an integer. The simplicity of the Griess algebra can be expected to
further restrict these possibilities. A full description of these
``exceptional Griess algebras'' and their possible realization in terms of
VOAs would be of obvious interest~\cite{T}.

Even stronger constraints on $\hat{V}^{(2)}$ follow if we assume that all $%
\mathrm{Aut}(V)$ invariants of $V^{(n)}$, $V^{(n)}/\mathrm{Aut}\,(V)$, are
Virasoro descendents of the vacuum for $n\leq 6$:$\,$

\begin{proposition}
Suppose that $V^{(n)}/\mathrm{Aut}\,(V)\subset V^{(n)}(C,0)$ for all $n\leq 6$. 
Then, $\hat{V}^{(2)}$ is an irreducible representation of $\mathrm{Aut}\,(V)$.
\end{proposition}

This concurs with the irreducible representation structure discussed in
Section 4.4 with $d_{2}(C=24)=196883$ for the Monster, $d_{2}(C=23\frac{1}{2}%
)=96255$ for the Baby Monster and $d_{2}(C=8)=155$ for $O_{10}^{+}(2)$.

The higher Casimirs lead to further constraints: if $\mu^{(n)}\in
V^{(n)}(C,0)$ for $n\leq 8$ or $n\leq 10$ then $C=24$ with $d_{2}=196883$.
However, $\mu^{(n)}\in V^{(n)}(C,0)$ for $n\leq 12$ is impossible.
This concurs with the fact that the Moonshine module contains a Monster
invariant primary vector of level~12.

We can also consider the Casimirs $\nu ^{(n)}=u_{3-n}^{\alpha }u_{\alpha
}\in V^{(n)}$ constructed from a basis $\{u_{\alpha }\}$ of primary vectors
of level $3$ 
\begin{gather*}
\hat{V}^{(3)}=V^{(3)}-L_{-1}V^{(2)},
\end{gather*}
of dimension $d_{3}=\dim \hat{V}^{(3)}>0$. We f\/ind

\begin{proposition}
Suppose that $\nu ^{(n)}\in V^{(n)}(C,0)$ for $n\leq 10$. Then $%
d_{3}(C)=p(C)/q(C)$ with 
\begin{gather*}
p(C) =5C\left( 5C+22\right) \left( 2C-1\right) \left( 7C+68\right) \left(
3C+46\right) \left( 5C+3\right) \left( 11C+232\right),  \\
q(C) =75{C}^{6}-9945{C}^{5}+472404{C}^{4}-9055068{C}^{3} \\
\phantom{q(C)=}{}+39649632{C}^{2}+438468672 C+2976768. 
\end{gather*}
\end{proposition}

All the factors of $\det M^{(10)}$ appear in $p(C)\,$explaining the remarks
concerning prime divisors in Table~2 where $d_{3}(24)=21296876=
2^{2}\cdot 31\cdot 41\cdot 59\cdot 71$ as obtains for the Moonshine Module. Finally, we also f\/ind

\begin{proposition}
Suppose that $V^{(n)}/\mathrm{Aut}\,(V)\subset  V^{(n)}(C,0)$ for all $n\leq
10$. Then $\hat{V}^{(3)}$ is an irreducible representation of $\mathrm{Aut}\,(V)$. 
\end{proposition}

This agrees with the Monster irreducible decomposition of $V^{(3)}$ for the
Moonshine Module where $\dim V^{(3)}=21493760=1+\dim \hat{V}^{(2)}+\dim \hat{V}^{(3)}=1+196883+21296876$.

\subsection*{Acknowledgments} 

The author thanks A.~Matsuo and G.~Mason for very
useful discussions and H.~Maruoka, A.~Matsuo and H.~Shimakura for generously
making Ref.~\cite{MMS} available.

\pdfbookmark[1]{References}{ref}
\LastPageEnding

\end{document}